\documentclass[12pt]{amsart}

\usepackage[cp1251]{inputenc} 
\usepackage[russian]{babel}
\usepackage{amsthm,amssymb,amsmath,amsfonts,mathrsfs}
\usepackage{mathbbol}
\theoremstyle{definition}

\theoremstyle{remark}

\usepackage{bbm}
\begin{document}
\normalsize
\begin{center}
\bf  Belyi pairs in the critical filtrations of Hurwitz spaces\footnote{based on the talk, given at the Workshop on Grothendieck-Teichm\"uller Theories 
(Chern Institute of Mathematics, Nankai University, Tianjin, China,
July 25 -- 29, 2016}\\
\end{center}
\begin{center}
\bf George SHABAT\\
 \vspace{.7cm}
\end{center}
\small
\bf0. Introduction\dotfill1\\
1. Hurwitz spaces\dotfill3\\
...1.0. Notations and conventions\dotfill3\\
...1.1. Definitions\dotfill4\\
...1.2. Main construction\dotfill5\\
...1.3. Critical filtration\dotfill8\\
2. Fried pairs and families\dotfill12\\
...2.0. Fried pairs\dotfill12\\
...2.1. Fried  families\dotfill13\\
...2.2. Structures on the base\dotfill13\\
...2.3. Dimension of the base\dotfill14\\
...2.4. Undesirable cases\dotfill15\\
...2.5. Belyi function on the base of a Fried family\dotfill15\\
3. Examples\dotfill16\\
...3.0.  Critical filtrations, genus 0\dotfill16\\
...3.1. Genus 1, degree 3\dotfill18\\
...3.2. Genus 2, degree 5\dotfill20\\
...3.3. Dreaming about applications\dotfill21\\
4. Conclusion\dotfill22\\
 References\dotfill 23
\rm
\normalsize
\vspace{3mm}
\begin{center}
\bf0. Introduction
\end{center}
\vspace{1mm}
This paper has a Chinese origin, and during the Tianjin conference the author (as, perhaps,  the other European and American participants) experienced an unusual feeling of being  surrounded by characters \it meaning \rm something  -- unlike the letters of latin or cyrillic alphabets. It was exciting to know\footnote{see, e.g., $\bf{[Xiao2009] }$} that such words as \it 
 year, 
 home, 
 sky (heaven),
 good,
 again, ...
\rm are representable by small beautiful pictures comprehensible to everybody (except a negligible group of foreigners).\\
\\
Do we have anything like that in mathematics? Can the topological objects be recognizable as the images of, say, arithmetical or algebraic objects?\\
\\
The affirmative answers, provided by 0- and 1-dimensional  cell complexes -- like
\begin{picture}(40,30)
\put(20,10){\circle*{3}}
\put(27,4){\circle*{3}}
\put(5,1){\circle*{3}}
\put(8,13){\circle*{3}}
\put(22,2){\circle*{3}}
\put(30,10){\circle*{3}}
\end{picture} 
, representing natural numbers, or
\begin{picture}(40,30)
\put(0,3){\line(1,0){40}}
\put(20,3){\circle*{3}}
\put(27,3){\circle*{3}}
\put(0,3){\circle*{3}}
\put(8,3){\circle*{3}}
\put(40,3){\circle*{3}}
\end{picture} 
and
\begin{picture}(40,30)
\put(20,10){\circle{30}}
\put(22,25){\circle*{3}}
\put(35,4){\circle*{3}}
\put(7,1){\circle*{3}}
\put(8,21){\circle*{3}}
\put(22,-6){\circle*{3}}
\put(31,22){\circle*{3}}
\end{picture} 
 representing linear and cyclic orders respectively, are not quite persuasive; perhaps, finite sets and graphs are too flexible to keep information about more complicated structures.\\
\\
However, things become more interesting in real dimension 2; as soon as we replace the abstract graphs by the \it embedded \rm  ones (into the closed oriented surfaces, in such a way that the complement should be homeomorphic to a disjoint union of open discs), we arrive at  the Grothendieck's theory of \it dessins d'enfants \rm -- see the original ideas in $\bf{[Grothendieck1984] }$, their mathematical treatment in $\bf{[ShabatVoevodsky1990] }$ and an elementary introduction in $\bf{[LandoZvonkin2004] }$. This theory is based on the equivalence between the category of purely algebraic objects, \it Belyi pairs\rm, (they occupy the considerable part of the present paper) and the above-mentioned category of  \it graphs on surfaces\rm; the latter may be considered as a category of \it images \rm of the objects of the former. One can daresay that  dessins d'enfants \rm or simply \it dessins\rm) serve like \it hieroglyphs\footnote{according to one of the vocabularies, a hieroglyph is \it any figure, character, or mark having or supposed to have a mysterious or enigmatic significance...} \rm  of arithmetic curves; see the discussion in $\bf{[Shabat2008] }$. \\
\\
In any language, even an artificial one aimed at the description of a certain class of objects, the \it grammar \rm  is at least as important as the \it alphabet\rm. One of the most fundamental structure on the set of (equivalence classes of) dessins is defined by the action of the absolute Galois group, see $\bf{[Grothendieck1984] }$;  the present paper is devoted to an hierarchy of other, ''grammatical-like'' structures, that realize certain fundamental relations between dessins.\\
\\
The paper is 
organized \normalsize as follows. Section 1 contains some foundational material, introducing  the versions of Hurwitz spaces needed for the further constructions, the main of which is the \it critical filtration. \rm  Section 2 is devoted to the 1-dimensional strata of this filtration. Section 3 contains several simplest examples of the introduced objects.\\
\\
The central concept of the paper is \it the  Belyi function on the base of a Fried family\rm, defined in the end of section 2. According to Grothendieck theory, the existence of such objects imply the possibility of encoding the main structures of Fried families graphically -- this is related to the above-mentioned \it grammar \rm of the language of dessins d'enfants. However, the paper contains only algebraic counter-parts of these structures -- representing them pictorially needs the elaboration of a special techniques of \it colored \rm drawing; it is (hopefully) postponed until further publications.\\
\\
The author is  indebted \normalsize to Mike Fried and Dima Oganesyan for useful discussions. Some critical remarks of the referee on the preliminary version of the paper resulted in the clarification of certain points.\\
\\
The  paper is supported in part by the Simons foundation. 
\\
\begin{center}
\bf1. Hurwitz spaces
\end{center}
Following $\bf{[Clebsch1872] }$ and $\bf{[Hurwitz1891] }$,  we are going to consider pairs $(X,f)$ of curves and rational functions on them; we are going to use the notation $f=x$ for the general cases (having in mind a curve $X$ with an affine model, defined by the polynomial equation $F(x,y)=0$) and $f=\beta$, $f=\varphi$ in the special cases defined below, where the letters $\beta$ and $\varphi$ are used in honor of Belyi and Fried. Though   several versions of Hurwitz spaces are considered in a great amount of papers, none of them fits our needs completely. E.g., in $\bf{[RomagnyWewers2006] }$ only the general enough functions are considered (as in the classical papers), while in $\bf{[Fried1977] }$ only the functions that are factorizations by the actions of finite groups are considered. The compactifications of classical Hurwitz spaces were studied, e.g., in, $\bf{[Deopurkar2014] }$, but their the generic functions on singular curves were added and not non-generic functions on smooth curves, the ones  we need.\\
\\
Therefore we suggest a brief self-consistent introduction to the theory of the extended Hurwitz spaces.\\
\\
\bf1.0. Notations and conventions. \rm
Let $\Bbbk$ be an algebraically closed field; by a \it curve \rm we mean a smooth irreducible complete curve over $\Bbbk$.\\
\\
Let $d\in\mathbb{N}\setminus\{0\}$ and $g\in\mathbb{N}$ be natural numbers. Denote
 $\mathcal{M}_g(\Bbbk)$ the \it moduli space \rm of curves of genus $g$. We use the standard notation
  $\mathcal{M}_{g,d}(\Bbbk)$ for the set of isomorphic classes of curves of genus $g$ with $d$ \underline{ordered} marked points and less standard 
   $\mathcal{M}_{g,[d]}(\Bbbk)=\mathcal{M}_{g,d}(\Bbbk)/\mathsf{S}_d$ for that of $d$ \underline{unordered} marked points.\\
  \\
  The field $\Bbbk$ will be fixed and the reference to it will often be omitted.\\
  \\
  For a curve $X$ of genus $g$ denote $[X]\in\mathcal{M}_g(\Bbbk)$ the corresponding point in the moduli space; 
 $\Bbbk(X)$ is the field of rational functions on $X$.\\
 \\
 We identify a function $x\in\Bbbk(X)$ with a mapping $x:X\to\mathbf{P}_1(\Bbbk)$. For $x\notin\Bbbk$ we define its \it degree \rm as the degree of the  field extension
$$
\deg x:=[\Bbbk(X):\Bbbk(x)].
$$
When this extension is \it separable \rm  (which is usually assumed) the degree of $x$ coincides with
the cardinality of pre-images\footnote{we use the notation $f^{-1\circ}$ for the \it compositional \rm inverse; thus, $\tan^{-1\circ}=\arctan$, while $\tan^{-1}=\cot$.}
$$
\deg x=\#x^{-1\circ}(P)
$$
for \underline{almost} all $P\in\mathbf{P}_1(\Bbbk)$.\\
\\
Denote the set of \it critical values \rm of a function $x\in\Bbbk(X)\setminus\Bbbk$
$$
\text{CritVal}(x):=\{P\in\mathbf{P}_1(\Bbbk)\mid\#x^{-1\circ}(P)<\deg x\}.
$$
 We use a bit nonstandard  notation  $\Bbbk(X)_d$ for the set of rational functions on $X$ of degree $d$.\\
\\
For a curve $X$ the group of its divisors is $\text{Div}(X)$ and for a function 
$f\in\Bbbk(X)\setminus\{0\}$ its divisor is $\text{div}(f)$. For a divisor $D\in\text{Div}(X)$ we'll use  the $\Bbbk$-vector space $\text{L}(D):=\{f\in\Bbbk(X)\mid\text{div}(f)+D\ge0\}$; according to the Riemann-Roch theorem, $\dim_\Bbbk\text{L}(D)=\deg D+1-g$ as soon as $\deg D> 2g-2$.\\
\\
\bf1.1. Definitions. \rm We consider two versions of Hurwitz spaces. Set-theoretically both consist of the equivalence classes of  pairs $(X,x)$ of a curve and a function on it, but modulo two different equivalence relations:
$$
\mathcal{HUR}_{d,g}^\nabla(\Bbbk):=\frac{\{(X,x)\mid[X]\in
\mathcal{M}_g(\Bbbk),x\in\Bbbk(X)_d\}}{\approx_\nabla}
$$
and
$$
\mathcal{HUR}_{d,g}^\square:=\frac{\{(X,x)\mid[X]\in
\mathcal{M}_g(\Bbbk),x\in\Bbbk(X)_d\}}{\approx_\square},
$$
where $(X,x)\approx_\nabla(X',x')$ iff there exists a commutative diagram
\begin{center}
\begin{picture}(110,60)
\put(28,50){$X$}
\put(53,-5){$\mathbf{P}_1(\Bbbk)$}
\put(35,45){\vector(1,-2){20}}
\put(39,54){\vector(1,0){52}}
\put(92,50){$X'$}
\put(95,45){\vector(-1,-2){20}}
\put(58,55){$\simeq$}
\put(35,22){$x$}
\put(88,22){$x'$}
\end{picture}
\end{center}
and $(X,x)\approx_\square(X',x')$ iff there exists a commutative diagram
\begin{center}
\begin{picture}(110,60)
\put(28,50){$X$}
\put(23,-5){$\mathbf{P}_1(\Bbbk)$}
\put(83,-5){$\mathbf{P}_1(\Bbbk)$}
\put(35,45){\vector(0,-1){40}}
\put(39,54){\vector(1,0){52}}
\put(55,-2){\vector(1,0){25}}
\put(92,50){$X'$}
\put(95,45){\vector(0,-1){40}}
\put(58,55){$\simeq$}
\put(25,22){$x$}
\put(98,22){$x'$}
\put(64,1){$T$}
\end{picture}
\\
\vspace{.5cm}
\end{center}
with some $T\in\mathsf{PSL}_2(\Bbbk)$, acting on  $\mathbf{P}_1(\Bbbk)$ by a fractional-linear transformation.\\
\\
Obviously, the group $\mathsf{PSL}_2(\Bbbk)$ of fractional-linear transformations acts on $\mathcal{HUR}_{d,g}^\nabla$ by post-composing:
$$
T\cdot[X,x]:=[X,T\circ x],
$$
and set-theoretically the space $\mathcal{HUR}_{d,g}^\nabla$ is \it almost \rm
 a $\mathsf{PSL}_2(\Bbbk)$-principal fiber bundle
$$
\mathcal{HUR}_{d,g}^\nabla
\stackrel{/\mathsf{PSL}_2}
\longrightarrow\mathcal{HUR}_{d,g}^\square;
$$
\it almost \rm because over curves with non-trivial automorphisms the fibers are further factorized over non-trivial finite groups.\\
\\
\bf1.2. Main construction. \rm From now on, considering the pair $(X,x)$, we assume that $x$ is \it separable. \rm The sufficient condition for it is $\deg(x)\notin\text{char}(\Bbbk)\mathbb{N}$.\\
\\
 First of all, we restrict our attention by the pairs $(d,g)$, for which the forgetful map\\
$$
\mathcal{HUR}_{d,g}^\bigstar\longrightarrow
\mathcal{M}_g:
(X,x)\mapsto X
$$
is surjective. According to the Brill-Noether theory (see $\bf{[GriffHarr1978]}$, p. 261 ), it happens if
$$
d\ge\lfloor \frac{g+3}{2}\rfloor; 
$$
 however, we are going to impose the stronger condition. \\
 \\
 Consider the Zariski-open subset 
 $$
\mathcal{HUR}_{d,g}^{\nabla,\text{ simple poles}}\subset
\mathcal{HUR}_{d,g}^{\nabla},
$$
consisting of such pairs $(X,x)$, that $x$ has only simple poles.\\
 \\
Then consider the map
$$
\text{poles}_{d,g}:\mathcal{HUR}_{d,g}^{\nabla,\text{ simple poles}}\longrightarrow\mathcal{M}_{g,[d]}:
(X,x)\mapsto(X,\text{poles}(x)),
$$
the fibers of which set-theoretically are
$$
\text{poles}_{d,g}^{-1\circ}(X,\{P_1,\dots,P_d\})=\text{L}(P_1+\dots P_d).
$$
As it was mentioned,  these spaces are nonzero for $
d\ge\lfloor \frac{g+3}{2}\rfloor
$ and are of equal  (i.e., independent  of  the set $\{P_1,\dots P_d\}$) dimension $d-g+1$ for $d>2g-2$.\\
\\
\bf Proposition. \it For $d>2g-2$ the set $\mathcal{HUR}_{d,g}^{\nabla,\text{ simple poles}}$ admits the structure of quasi-projective manifold, with respect to which the map $\mathrm{poles}_{d,g}$ 
is a projection of vector bundle over the moduli space.\\
\\
\bf Sketch of the proof. \rm The problem is to locally trivialize the union of vector spaces
 $\text{L}(P_1+\dots P_d)$ over some (small enough) Zariski-open subsets of $\mathcal{M}_{g,[d]}$. It is solved constructively in many papers, e.g., in $\bf{[Hess2002] }$.\\
 \\
 A non-constructive proof can look as follows. For an arbitrary 
 $(X,\{P_1,\dots,P_d\})\in\mathcal{M}_{g,[d]}$ choose local parameters  $x_1,\dots,x_d$ around $P_1,\dots, P_d$. If we consider the tricanonical embedding of \it all \rm the curves of genus $g$ into 
 $\mathbf{P}_{5g-6}$, these local parameters can be chosen as ratios of appropriate linear combinations of homogeneous coordinates in this ambient projective space. Then the differential of a local coordinate $x_i$ remains nonzero in some Zariski-open neighborhood $U\subset\mathcal{M}_{g,[d]}$ of the chosen point
 $(X,\{P_1,\dots,P_d\})\in\mathcal{M}_{g,[d]}$, and, understanding $x_i's$ as the above expressions, we can construct local coordinates $P_i'\mapsto x_i-x_i(P_i')$ for all $(X',\{P_1',\dots,P_d'\})\in U$. \\
 \\
 Now the task of constructing the global rational functions with simple poles over $U$ can be understood as the solution of \it Cousin problem \rm with principal parts $\sum_{i=1}^d\frac{a_i}{x_i}$. Due to the exact sequences of sheafs
 $$
 0\longrightarrow\mathcal{O}_{X'}\longrightarrow\Bbbk(X')\longrightarrow
 \frac{\Bbbk(X')}{\mathcal{O}_{X'}}\longrightarrow0
 $$
 all over $U$, the above Cousin problem is solvable under $g$ independent linear conditions on the tuples $(a_1,\dots,a_d)$, resulting from the surjections
 $$
\text{H}^0\Big(X', \frac{\Bbbk(X')}{\mathcal{O}_{X'}}\Big) \longrightarrow
\text{H}^1(X', \mathcal{O}_{X'}).
 $$
However, the spaces $\text{H}^1(X', \mathcal{O}_{X'})$ are Serre dual to the spaces of regular differentials $\Omega^1[X']$, the union of which over $U$ constitute an open part of the well-known vector bundle, since it is a pullback of the \it Hodge bundle \rm under the forgetful morphism 
$\mathcal{M}_{g,[d]}\to\mathcal{M}_g$.\\
\\
Hence the Cousin problem under discussion reduces to the linear algebra problem that has just been locally trivialized over $\mathcal{M}_{g,[d]}$. $\blacksquare$\\
\\
The construction of the  proposition does not cover the Zariski-closed subset of $\mathcal{HUR}_{d,g}^\nabla$, consisting of the pairs $(X,x)$, for which $x$ has poles of order $\ge2$. However, the image of its complement $\mathcal{HUR}_{d,g}^{\nabla,\text{ simple poles}}$ in  $\mathcal{HUR}_{d,g}^{\square}$ is obviously the whole $\mathcal{HUR}_{d,g}^{\square}$: a generic fractional-linear transformation of any non-constant function has only simple zeros.\\
\\
Thus our quasi-projective structure can be pushed down from $\mathcal{HUR}_{d,g}^{\nabla,\text{ simple poles}}$ to the whole  $\mathcal{HUR}_{d,g}^{\square}$; then we can lift it back to the whole
 $\mathcal{HUR}_{d,g}^{\nabla}$.\\
 \\
 Our considerations are obviously consistent with the traditional ones, where $(X,x)$'s are considered for 
 $x$'s of Morse type; we just extend the well-known constructions to arbitrary $x$'s.\\
 \\
 The resulting formulation\footnote{from which we exclude the trivial cases $(d,g)=(1,0)$ and $(d,g)=(2,0)$} is as follows.\\
 \\
 \bf Theorem. \it Let d and g be natural numbers, satisfying $d\ge3$ for $g=0$ and  $d\ge2g+1$ for positive g's. Then the  sets $\mathcal{HUR}_{d,g}^\nabla$ and $\mathcal{HUR}_{d,g}^\square$ have the natural structures of quasi-projective manifolds of dimensions $2(d+g)-2$ and $2(d+g)-5$. The action of $\mathsf{PSL}_2$ on $\mathcal{HUR}_{d,g}^\nabla$, defined by 
 $T\cdot[X,x]_{\approx_\nabla}=[X,T\circ x]_{\approx_\nabla}$, is regular; its  stationary groups are finite -- with the exception of the pairs, isomorphic to $(\mathbf{P}_1,z\mapsto z^n)$. The factorization over this action coincides with the morphism
 $\mathcal{HUR}_{d,g}^\nabla\to\mathcal{HUR}_{d,g}^\square:[X,T\circ x]_{\approx_\nabla}
 \mapsto[X,T\circ x]_{\approx_\square}$. \\
 \\
 \bf Sketch of proof. \rm Most of statements  follow from definitions and the above constructions.  The stationary group of $\mathsf{PSL}_2$, fixing $[X,x]_{\approx_\nabla}$, is
 $$
 \mathrm{Aut}(X,x):=\{\alpha\in\mathrm{Aut}X\mid\exists T\in\mathsf{PSL}_2\text{ such that }
 T\circ x\circ\alpha\equiv x\}.
 $$ 
For $g\ge2$ this group is finite just because $\mathrm{Aut}X$ is finite. For $g=0$ and $g=1$ the finiteness of the group follows from the consideration of critical values of $x$. Their number is at least 3: for $g=0$ due to the restriction concerning  $z\mapsto z^n$ and for $g=1$ by obvious general reasons.\\
\\
As for dimensions, we have realized $\mathcal{HUR}_{d,g}^\nabla$ as a total space of vector bundle,
so
$$
\dim\mathcal{HUR}_{d,g}^\nabla=
\dim{\mathcal{M}_{g,d}}+
\dim\mathrm{L}_X(P_1+\dots+P_d)=_{\text{Riemann-Roch}}
$$
$$
=(3g-3+d)+d-g+1=2(d+g)-2.
$$
Due to the $\mathsf{PSL}_2$ statement,
$$
\dim\mathcal{HUR}_{d,g}^\square=\dim\mathcal{HUR}_{d,g}^\nabla-3=2(d+g)-5. \blacksquare
$$
\\
In the case $\Bbbk=\mathbb{C}$ these dimensions of Hurwitz spaces are well known and easily obtainable
by the transcendental theory of ramified coverings.\\ 
\begin{center}
\bf1.3. Critical filtration
\end{center}
From now on we are going to consider both types of Hurwitz spaces introduced and use
$$
\bigstar\in\{\nabla,\square\}.
$$
A  special object of our considerations
is the \it critical filtration \rm
$$
\mathcal{HUR}_{d,g}^\bigstar=\mathcal{HUR}_{d,g;2(d+g)-2}^\bigstar
\supset\mathcal{HUR}_{d,g;2(d+g)-3}^\bigstar\supset\dots,
$$
where
$$
\mathcal{HUR}_{d,g;\mathbf{b}}^\bigstar:=\{[(X,x)]_\bigstar\mid
\#\mathrm{CritVal}(x)\le\mathbf{b}\}.
$$
We want the critical filtrations to be \underline{as long as possible}; hence the question arises about the smallest possible value of $\mathbf{b}$ in the last (non-empty) term 
$\mathcal{HUR}_{d,g;\mathbf{b}}^\bigstar$ of the critical filtration. In other words: 
\begin{center}
\bf how small can be \rm 
$\#\text{CritVal}(x)$ for $x\in\Bbbk(X)\setminus\Bbbk$?
\end{center}
The answer is simple.\\
\\
$\#\text{CritVal}(x)=0$ iff $(X,x)\simeq_\square(\mathbf{P}_1,z\mapsto z)$;\\
$\#\text{CritVal}(x)\ne1$ never;\\
$\#\text{CritVal}(x)=2$ iff $(X,x)\simeq_\square(\mathbf{P}_1,z\mapsto z^d)$ for some 
$d\in\mathbb{N}$ -- see Lemma inside the proposition below;\\
$\#\text{CritVal}(x)=3$ in lots of cases; such $(X,x)$'s are called \it Belyi pairs\rm.\\
\\
If $(X,x)$ is a Belyi pair, then  $x$ is called a \it Belyi function \rm on a curve $X$. According to the famous Belyi theorem (see $\bf{[Belyi1980] }$, $\bf{[Belyi2002] }$), there exists a Belyi function on \underline{any} curve over $\overline{\mathbb{Q}}$.\\
\\
Hence we have to define the pairs
 $(d,g)$, for which the critical filtration is \it complete\rm, i.e.
$$
\mathcal{HUR}_{d,g}^\bigstar=\mathcal{HUR}_{d,g;2(d+g)-2}^\bigstar
\supset\mathcal{HUR}_{d,g;2(d+g)-3}^\bigstar\supset\dots
\supset\mathcal{HUR}^\bigstar_{d,g;3}\ne\varnothing;
$$
The answer turns out to be familiar.\\
\\
\bf Proposition. \it If the stratum $\mathcal{HUR}^\bigstar_{d,g;3}$ is non-empty, then 
$d\ge2g+1$. Conversely, if $d\ge2g+1$ and $\mathrm{char}(\Bbbk)$ is zero or large enough, then
the stratum $\mathcal{HUR}^\bigstar_{d,g;3}$ is non-empty.\\
\\
 \bf Proof. \rm First suppose that the stratum $\mathcal{HUR}^\bigstar_{d,g;3}$ is non-empty; it means that there exists a Belyi pair\footnote{we use the notation for $\beta$ standard since $\bf{[ShabatVoevodsky1990]}$}
$(X,\beta)\in\mathcal{HUR}^\nabla_{d,g}$. Applying to $\beta$ an appropriate fractional-linear transformation, we can assume that $\text{CritVal}(\beta)\subseteq\{0,1,\infty\}$.\\
\\
Consider the divisors
$$
\text{div}(\beta)=\sum_{i=1}^{\alpha_0} a_iA_i-\sum_{k=1}^{\alpha_\infty} c_kC_k\eqno{\bf{(1.3.1)}}
$$
and, if 1 is the third critical value of $\beta$,
$$
\text{div}(\beta-1)=\sum_{j=1}^{\alpha_1} b_jB_j-\sum_{k=1}^{\alpha_\infty} c_kC_k,\eqno{\bf{(1.3.2)}}
$$
where $\alpha_0,\alpha_1,{\alpha_\infty}\in\mathbb{N}$. We assume that $\beta$ is non-constant, hence 
$$
\alpha_0>0,{\alpha_\infty}>0, \eqno{\bf{(1.3.3)}}
$$
 and $A_1,\dots,A_{\alpha_0},B_1,\dots,B_{\alpha_1},C_1,\dots,C_{\alpha_\infty}\in X$ are \underline{different} points, 
 while $a_1,\dots,a_{\alpha_0},b_1,\dots,b_{\alpha_1},c_1,\dots,c_{\alpha_\infty}\in\mathbb{N}_{>0}$. By definition,
 $$
 \sum_{i=1}^{\alpha_0} a_i=\sum_{k=1}^{\alpha_\infty} c_k=d. \eqno{\bf{(1.3.4)}}
 $$
 For the time being we admit the possibility ${\alpha_1}=0$, and the case $\#\text{CritVal}=2$ has to be considered separately. \\
 \\
\bf Lemma. \it  If $n=0$, then $(X,x)\simeq_\square(\mathbf{P}_1,z\mapsto z^d)$.\\
\\
\bf Proof. \rm Since 0 and $\infty$ are supposed to be the only critical values of $\beta$, the zeros and poles of its differential are contained in the set of its zeros and poles and hence\footnote{note the fonts: $d$ is \underline{a} natural number while d is \underline{the} exterior differential}
$$
\text{div}(\text{d}\beta)=\sum_{i=1}^{\alpha_0} (a_i-1)A_i
-\sum_{k=1}^{\alpha_\infty} (c_k+1)C_k.\eqno{\bf{(1.3.5)}}
$$
Then
$$
2g-2=\deg\text{div}(\text{d}\beta)=\sum_{i=1}^{\alpha_0} (a_i-1)
-\sum_{k=1}^{\alpha_\infty} (c_k+1)=
$$
$$
=\big(\sum_{i=1}^{\alpha_0} a_i\big)-\alpha_0
-\big(\sum_{k=1}^{\alpha_\infty} c_k\big)-{\alpha_\infty}=_{\mathbf{(1.3.4)}}-\alpha_0-{\alpha_\infty},
$$
so  \bf(1.3.3) \rm implies $2g-2\le-2$. Hence $g=0$, and a rational function with only two critical values is $\mathsf{PSL}_2$-equivalent to $z\mapsto z^d. \blacksquare$\\
\\
\bf Remark. \rm Over $\Bbbk=\mathbb{C}$ the above lemma is immediate, since the assumption 
$\#\text{CritVal}(\beta)=2$ implies that the fundamental group 
$\pi_1(\mathbf{P}_1(\mathbb{C})\setminus\text{CritVal}(\beta))$ is $\mathbb{Z}$ and therefore all the non-ramified coverings of $\mathbf{P}_1(\mathbb{C})\setminus\text{CritVal}(\beta)$ are of the desired type.\\
\\
In the general case the inequality
$$
n>0\eqno{\bf{(1.3.6)}}
$$
is added to \bf(1.3.3)\rm, and by the same reason as in the Lemma (the zeros and poles of the differential of a Belyi function are among the pre-images of its critical values) we have\\
$$
\text{div}(\text{d}\beta)=\sum_{i=1}^{\alpha_0} (a_i-1)A_i+\sum_{j=1}^{\alpha_1}( b_j-1)B_j
-\sum_{k=1}^{\alpha_\infty} (c_k+1)C_k. \eqno{\bf{(1.3.7)}}
$$
Hence
$$
2g-2=\deg\text{div}(\text{d}\beta)=\sum_{i=1}^\alpha (a_i-1)+\sum_{j=1}^n( b_j-1)
-\sum_{k=1}^\gamma (c_k+1)=
$$
$$
=\big(\sum_{i=1}^{\alpha_0} a_i\big)-\alpha_0+\big(\sum_{j=1}^{\alpha_1} b_j\big)-{\alpha_1}
-\big(\sum_{k=1}^{\alpha_\infty} c_k\big)-{\alpha_\infty}=
$$
$$
=d-\alpha_0-{\alpha_1}-{\alpha_\infty}.
$$
So
$$
d=2g-2+\alpha_0+{\alpha_1}+{\alpha_\infty}\ge_{\mathbf{(1.3.3),(1.3.6)}}2g-2+3=2g+1,
$$
and half of the proposition is proved.\\
\\
For the other half, consider first the curve
$$
y^2=1-x^{2g+1};
$$
it has genus $g$, if $\text{char}(\Bbbk)\ne2$ and  $2g+1\notin\text{char}(\Bbbk)\mathbb{N}$.\\
\\
Consider the function
$$
\beta=y
$$
on this curve; obviously, under the same characteristic assumptions $\deg\beta=2g+1$.\\
\\
Calculating
$$
\mathrm{d}\beta=\mathrm{d}y=-\frac{2g+1}{2}x^{2g}\frac{\mathrm{d}x}{y}
$$
and remembering that the differentials $x^{k}\frac{\mathrm{d}x}{y}$ for $0\le k\le g-1$ are regular on the curve, we see that  $\mathrm{d}\beta=0\Longrightarrow x=0\Longrightarrow y=\pm1$, so $\beta$ is a Belyi function with critical values $\{\pm1,\infty\}$.\\
\\
We have presented a sequence of Belyi pairs, covering the \it extreme \rm values $d=2g+1$ of our proposition; it remains to show that the degrees  of Belyi functions can be increased arbitrarily within a fixed genus.\\
\\
The corresponding statement is totally obvious over $\Bbbk=\mathbb{C}$ because of the correspondence between the Belyi pairs and the \it dessins d'enfants \rm -- see $\bf{[LandoZvonkin2004] }$; the needed statement reduces to the possibility of attaching edges to embedded graphs. \\
\\
Due to the \it Lefshetz principle \rm (see, e.g., $\bf{[ Bouscaren1998] }$) the proposition holds over any $\Bbbk$ of characteristic 0, and due to the \it transfer principle \rm (see $\bf{[Robinson1963] }$) also
over $\Bbbk$'s of large enough characteristic.\\
\\
Some restrictions on the characteristic are necessary; e.g., there is the only Belyi pair with $(d,g)=(3,1)$ (see below), and it does not admit reduction modulo 3. $\blacksquare$\\
\\
Now we have at least two reasons to stick to the inequality 
$$
\boxed{ d\ge2g+1}
$$
 Under this assumption we introduce the finite sets
$$
\mathrm{BEL}_{g,d}(\Bbbk):=\mathcal{HUR}_{d,g;3}^\square(\Bbbk)
$$
and the system of curves, the name to which is given because of the paper $\bf{[Fried1990] }$
$$
\mathrm{FRIED}_{g,d}(\Bbbk):=\mathcal{HUR}_{d,g;4}^\square(\Bbbk).
$$
A tremendous number of papers, including the most recent ones, is devoted to the calculation of
 \normalsize
 $\#\mathrm{BEL}_{g,d}(\mathbb{C})$. 
 On the contrary, the combinatorial structure of  \normalsize $\mathrm{FRIED}_{g,d}$ 
 does not seem to have been investigated enough\footnote{the author is indebted to the referee for pointing out 
 the paper $\bf{[BaileyFried2002]}$ containing some important information on the geometry of Fried curves. The concept of \it modular towers \rm should be related to the concepts of the present paper.}. The constructions of the present paper stress that 
 the latter constitutes the system of algebraic varieties 
 (in the Hurwitz or moduli spaces) 
 along which the points of the former are distributed according to some fundamental principles.
 \normalsize
 \\
 \begin{center}
 \bf2. Fried pairs and families
 \end{center}
 In this section we turn to the one-dimensional strata of the critical filtrations.\\
 \\
 \bf2.0. Fried pairs. \rm Recall that for a curve $X$ over $\Bbbk$, a natural number 
 $d\in\mathbb{N}\setminus\mathrm{char(\Bbbk)\mathbb{N}}$ and a function
 $\beta\in\Bbbk(X)_d$ if $[(X,\beta)]_\square\in\mathrm{BEL}_{g,d}$, then $(X,\beta)$ is called a 
 \it Belyi pair\rm;
 similarly, if $[(X,\varphi)]_\square\in\mathrm{FRIED}_{g,d}$, then $(X,\varphi)$ is called a 
 \it Fried pair\rm.\\
\\
\bf Examples of Belyi and Fried pairs\rm\\
\\
(a) If $X$ is a projective normal model of the \it generalized Fermat affine curve\rm, defined by the equation
$$
x^m+y^n=1,
$$
then
 $(X,x^m)=(X,1-y^n)$ is a Belyi pair. We met a particular case of it in the proof of proposition \bf1.3\rm.\\
\\
(b) For a \it Klein quartic \rm $X\subset\mathbf{P}_2(\Bbbk)$, defined by the  equation
$$
x^3y+y^3z+z^3x=0,
$$
the pair\footnote{recall that we identify rational functions an a curve  and (ramified) coverings of the projective line} $(X,X\to\frac{X}{\mathrm{Aut}X})$ is a Belyi one.\\
\\
(c)  Denote
$$
\ddot{\Bbbk}:=\Bbbk\setminus\{0,1\}.
$$ 
For $t\in\ddot{\Bbbk}$ and the elliptic curve $\mathbf{E}_t$, defined by the affine equation
$$
y^2=x(x-1)(x-t),
$$
the pair $(\mathbf{E}_t,x)$ is a Fried one. Here $(d,g)=(2,1)$, so the pair $(d,g)$ does not satisfy our main inequality $d\ge2g+1$. As a result, the fibers of this family do not degenerate to Belyi pairs. \\
\\
\bf2.1. Fried  families. \rm However, the concept of \underline{individual} Fried pair seems mathematically senseless -- these pairs should vary in the \underline{families}, as in the above example of elliptic curves.\\
\\
\bf Definition. \rm A \it Fried family over a field $\Bbbk$ \rm is a quadruple $(\mathcal{X},\mathbf{B},\pi,\Phi)$ of quasi-projective manifold 
$\mathcal{X}$ together with a proper flat morphism with smooth fibers $\pi:\mathcal{X}\longrightarrow \mathbf{B}$ onto a smooth quasi-projective irreducible manifold $\mathbf{B}$  and  a
morphism $\Phi:\mathcal{X}\longrightarrow\mathbf{P}_1(\Bbbk)$, called the \it Fried map\rm. Denoting for $B\in\mathbf{B}$ the fibers 
$$
\mathcal{X}_B:=\pi^{-1}(B),
$$
we assume that all its restrictions
to the fibers 
$$
\Phi|_{\mathcal{X}_B}:\mathcal{X}_B\longrightarrow\mathbf{P}_1(\Bbbk)
$$
are Fried functions.\\
\\
\bf Example. \rm In order to describe  the above family of elliptic curves $\mathbf{E}_t$ in the introduced terms, consider 
the  manifold $\mathcal{X}\subset\ddot{\Bbbk}\times\mathbf{P}_2(\Bbbk)$ over the base 
$\mathbf{B}:=\ddot{\Bbbk}$, defined by the equation
$$
y^2z=x(x-z)(x-tz).
$$
Its projection is
$$
\pi:\mathcal{X}\longrightarrow \mathbf{B}:(t,x:y:z)\mapsto t
$$
and the Fried map is defined by
$$
\Phi:\mathcal{X}\longrightarrow \mathbf{P}_1(\Bbbk):(t,x:y:z)\mapsto
\left\{
\begin{array}{l}
 ( x:z),\text{ if }(x:y:z)\ne(0:1:0)\\
  (1:0), \text{ if }(x:y:z)=(0:1:0).
\end{array}\right.
$$
\\
\bf2.2. Structures on the base. \rm
A Fried family over a base $\mathbf{B}$ defines  two structures on  $\mathbf{B}$:
one is related to to the isomorphic class of the fibers
$\mathcal{X}_B$'s and the other to the critical values of $\Phi|_{\mathcal{X}_B}$'s.\\
\\
The first one, depending only on $\pi$, is realized by a \it modular morphism\rm
$$
\mathbf{B}\longrightarrow\mathcal{M}_{g}(\Bbbk):B\mapsto[\mathcal{X}_B],
$$
and we postpone its analysis until the further publications. In the above example of family of elliptic curves $\mathbf{E}_t$ it is realized by the $\mathbf{j}$-invariant 
$$
t\mapsto256\,{\frac { \left( {t}^{2}-t+1 \right) ^{3}}{{t}^{2} \left( t-1 \right) ^{2}}}.
$$
The second is  defined only on the dense Zariski-open subset of $\mathbf{B}$ of points, over which the critical values do not collide, i.e. on the set
$$
\{B\in\mathbf{B}\mid \#\mathrm{CritVal}(\Phi|_{\mathcal{X}_B})=4\}.
$$
This rational map, denoted by $\beta_{\text{bas}}$ for the reasons that will be clarified soon\footnote{when we restrict our considerations by 1-dimensional bases}, is defined by
$$
\beta_{\text{bas}}:\mathbf{B}\dashrightarrow\mathcal{M}_{0,[4]}(\Bbbk):B\mapsto
\mathrm{CritVal}(\Phi|_{\mathcal{X}_B}).
$$
In order to replace the space of unordered quadruples $\mathcal{M}_{0,[4]}$ by that of the ordered ones $\mathcal{M}_{0,4}$ we have to lift this map to some cover of the original base
$$
\widehat{\beta_{\text{bas}}}:\widehat{\mathbf{B}}\dashrightarrow\mathcal{M}_{0,4}(\Bbbk)\simeq\Bbbk.
$$
The universal way to do it is to use the fiber product of $\beta_{\text{bas}}$ with the order-forgetting map $\mathcal{M}_{0,4}\to\mathcal{M}_{0,[4]}$, but the result will be reducible and some recipe for choosing an irreducible component will be needed.\\
\\
\bf2.3. Dimension of the base. \rm Sometimes it seems reasonable to consider multidimensional Fried families with  multidimensional bases. For example, the Weierstrass family with the affine models of fibers defined by
$$
\mathcal{X}_{g_2,g_3}:y^2=4x^3-g_2x-g_3,
$$
has the base $$\mathbf{B}=\{(g_2,g_3)\in\Bbbk\times\Bbbk\mid
{g_{{2}}}^{3}\ne27\,{g_{{3}}}^{2}\}.$$
It is a Fried family  with the Fried map $\Phi=x$ again.  \\
\\
The sets of isomorphic curves in this family are the orbits of the multiplicative group $\Bbbk^\times$, acting on the base by $\lambda\cdot(g_2,g_3):=(\lambda^3g_2,\lambda^2g_3)$; the fibers of the family over these orbits provide the isomorphic Fried pairs.\\
\\
This family is convenient for the study of elliptic curves, but for our purposes of classifying the bases of Fried families inside the Hurwitz or moduli spaces we would prefer its one-dimensional subfamilies with varying Fried pairs. Unfortunately, such subfamilies, containing all the curves \it exactly once\rm, do not exist; for example, over the line $g_3=0$ all the curves are isomorphic, while over the line $g_3=1$ we do not find the important curve $y^2=x^3-x$. This non-existence is a manifestation of the lack of the \it fine moduli space \rm of curves. 
 The best we can reach is to work with the families like $\{\mathbf{E}_t\}$ above, that contain all the curves, each one over the \it finite \rm  number of points of the base.\normalsize\\
\\
Generally, for any $(d,g)\in\dot{\mathbb{N}}\times\mathbb{N}$ consider any Fried family $\pi:\mathcal{X}\to\mathbf{B}$ with the fibers of genus $g$ and the Fried map of degree $d$ on the fibers. Then for every $B\in\mathbf{B}$   we have 
$[\mathcal{X}_B,\Phi|_{\mathcal{X}_B}]\in\mathcal{HUR}_{d,g;4}$, 
and the modular map  becomes the part of the commutative diagram 
\begin{center}
\begin{picture}(110,50)
\put(5,15){$\mathbf{B}$}
\put(65,44){$\mathcal{HUR}_{d,g;4}$}
\put(20,25){\vector(2,1){40}}
\put(20,20){\vector(1,0){44}}
\put(77,40){\vector(0,-1){16}}
\put(70,15){$\mathcal{M}_{g}$,}
\end{picture}
\end{center}
where the sloped arrow is $B\mapsto[\mathcal{X}_B,\Phi|_{\mathcal{X}_B}]$, and the vertical one is the forgetting map $[X,x]\mapsto[X]$. These sloped arrows constitute the \it etale topology \rm over the components of $\mathcal{HUR}_{d,g;4}$, that is (hopefully) the main object of the future theory.\\
\\
\bf2.4. Undesirable cases. \rm It may happen that the fibers of $\pi$ are isomorphic, while the Fried functions vary. This phenomenon  is, of course, inevitable in the case $g=0$. The example for $g>0$ is provided by the constant  ''family''
$y^2=(x^2-1)(x^2-2)$  with $\Phi_b=x^2-b^2$ and hence $\mathrm{CritVal(\Phi_b)}=\{\infty,-b^2,1-b^2,2-b^2\}$. In the rest of the paper  we avoid considering such families with fibers of positive genus.\\
\\
\bf2.5. Belyi function on the base of a Fried family. \rm 
The four critical values of Fried functions are unnumbered. Using cross-ratio 
$$
\langle a,b,c,d\rangle:=\frac{a-c}{b-c}\cdot\frac{b-d}{a-d}
$$
and numbering the critical values arbitrary
$$
\mathrm{CritVal}(\Phi_B)=:\{c_1,c_2,c_3,c_4\},
$$
we get the 6-valued ''function''
$$
B\mapsto t:=\langle c_1,c_2,c_3,c_4\rangle.
$$
Under the 24 re-numerations of critical values this quantity assumes only six values
$$
t,\frac{1}{t},1-t,\frac{1}{1-t},\frac{t}{t-1},\frac{t-1}{t},
$$
so  this quantity becomes the true function $\widehat{\mathbf{B}}\to\Bbbk$ only being lifted to the cover $\widehat{\mathbf{B}}\to\mathbf{B}$ of degree 1,2,3 or 6.\\
\\
It is easy to check that the above cross-ratios assume values $0,1,\infty$ iff the critical values collide, so the above partial function extends to a true one
$$
\widehat{\beta_{\mathrm{bas}}}:\widehat{\mathbf{B}}\longrightarrow\mathbf{P}_1(\Bbbk).
$$
We arrive at the main result of this section.\\
\\
\bf Theorem. \it For any Fried family $\pi:\mathcal{X}\to\mathbf{B}$ with the Fried function $\Phi$ the above-defined function $\widehat{\beta_{\mathrm{bas}}}$ on the appropriate cover of the base is a  Belyi function. Over any critical point $A\in\widehat{\mathbf{B}}$, i.e. a point, where $\widehat{\beta_{\mathrm{bas}}}(A)\in\{\infty,0,1\}$, the restriction of the Fried function $\beta_A:=\Phi|_{\mathcal{X}_A}$\it is a  Belyi function. \rm\\
\\
An alternative way of constructing the Belyi function on the base of a Fried family is to construct over each $B\in\mathbf{B}$ the elliptic curve $\mathbf{E}_B$, ramified over 
$\text{CritVal}(\Phi_B)\subset\mathbf{P_1(\Bbbk)}$ and then to use $\bf{j}$-invariants, defining
$$
\beta_{bas}(B):=\bf{j}(\bf{E}_B).
$$
In the remaining part of the paper we'll refer to this method as to \it elliptic trick\rm.\\
\\
Let us resume our constructions from the viewpoint   of dessins d'enfants theory. Having relaxed the condition \normalsize $\#\mathrm{CritVal}=3$ to $\#\mathrm{CritVal}=4$,  we still have Belyi pairs -- this time on the bases of Fried families -- and hence the opportunity to relate combinatorial topology, complex analysis and arithmetical geometry. The present paper is just the first step in realizing this possibility.\\
\begin{center}
\bf3. Examples
\end{center}
Usually the attempts to describe Fried families explicitly produce terribly cumbersome results. In several simplest cases below, however, the answers are relatively compact. The main results (they are boxed) contain Belyi functions on the bases of Fried families; their combinatorial and arithmetico-geometrical analysis will appear elsewhere.\\
\\
\bf3.0.  Critical filtrations, genus 0. \rm This case gives nothing for the study of moduli spaces by critical filtrations, but  provides some explicit examples of the components of strata of critical filtrations and their intersections.\\
\\
For all $d\in\mathbb{N}_{\ge3}$ there are isomorphisms 
$$
\mathcal{HUR}_{d,0}^\nabla(\Bbbk)\cong\frac{\Bbbk(x)_d}{\mathsf{PSL}_2(\Bbbk)}
$$
and
$$
\mathcal{HUR}_{d,0}^\square(\Bbbk)\cong\frac{\Bbbk(x)_d}
{\mathsf{PSL}_2(\Bbbk)\times\mathsf{PSL}_2(\Bbbk)}
$$
between  the Hurwitz spaces and the spaces of rational functions modulo groups of fractional-linear transformations; the group $\mathsf{PSL}_2(\Bbbk)$ acts on the space $\Bbbk(x)_d$ by $T\cdot R:=T\circ R$, and the group $\mathsf{PSL}_2(\Bbbk)\times\mathsf{PSL}_2(\Bbbk)$ -- by  $(T_1,T_2)\cdot R:=T_1\circ R\circ T_2^{-1\circ}$. So
$$
\dim\mathcal{HUR}_{d,0}^\nabla=2d-2
$$
and
$$
\dim\mathcal{HUR}_{d,0}^\square=2d-5.
$$
The geometry of critical filtration of these spaces can be described in terms of subvarieties of \it weighted projective spaces\rm; without going into details, describe these filtrations  briefly for $d=3$ and $d=4$.\\
\\
$\underline{d=3}$. Here $\dim\mathcal{{HUR}}_{3,0}^\square=1$, so the whole Hurwitz space is a Fried family, carrying the 2-term filtration 
$$
\mathcal{{HUR}}_{3,0;4}^\square\supset\mathcal{{HUR}}_{3,0;3}^\square.
$$
The functions can be normalized to the form
$$
R={\frac {{z}^{3}+{z}^{2}}{\lambda\,z+1}};
$$
The obviously special values are $\lambda=1$ (degree drops) and $\lambda=0$ (the poles collide). For  general values of $\lambda$ the critical values of $R$ are $R=\infty$, $R=0$ and the roots of the qudratic equation
$$
4\,{\lambda}^{3}{R}^{2}+ \left( {\lambda}^{2}+18\,\lambda-27 \right) R
+4,
$$
and these roots can not be parametrized by $\lambda$ rationally. \\
\\
The third one is $\lambda=9$, since with this value of parameter the function $R$ acquires a special critical point $z=-\frac{1}{3}$ with the critical value $R=-\frac{1}{27}$:
$$
{\frac {{z}^{3}+{z}^{2}}{9\,z+1}}+\frac{1}{27}=\frac{1}{27}\,{\frac { \left( 3\,z+1
 \right) ^{3}}{9\,z+1}}.
$$
Applying the \it elliptic trick\rm, we get the Belyi function
$$
\boxed{\beta_{\text{bas}}(\lambda)={\frac { \left( \lambda-3 \right) ^{3} \left( {\lambda}^{3}
-9\,{\lambda}^{2}+243\,\lambda-243 \right) ^{3}}{ \left( \lambda-1
 \right)  \left( \lambda-9 \right) ^{3}{\lambda}^{6}}}}
$$
Though very simple, this case is a bit enigmatic: it is mentioned in $\bf{[LandoZvonkin2004, \rm pp. 331-332)] }$ that the family under consideration is isomorphic to the family of \it polynomials \rm of degree 4, and that no algebro-geometric explanation of this phenomenon is known. Some explanations will be published in $\bf{[DremShab2017] }$.\\
\\
$\underline{d=4}$. Now $\dim\mathcal{{HUR}}_{4,0}^\square=3$, so we have the filtration 
$$
\mathcal{{HUR}}_{4,0;6}^\square\supset\mathcal{{HUR}}_{4,0;5}^\square
\supset{\mathcal{HUR}}_{4,0;4}^\square\supset\mathcal{{HUR}}_{4,0;3}^\square.
$$
The functions can be normalized to the form
$$
R={\frac {{z}^{4}+p{z}^{3}+q{z}^{2}}{r{z}^{2}+sz+1}}
$$
What survives at such $R$'s from the action of $\mathsf{PSL}_2(\Bbbk)\times\mathsf{PSL}_2(\Bbbk)$ on $\Bbbk(x)_4$ is the product of two small groups. The group $\Bbbk^\times$ acts upon these tuples by
$$
 \lambda\cdot (p,q,r,s,z,R):=(\lambda\,p,{\lambda}^{2}q
,{\frac {r}{{\lambda}^{2}}},{\frac {s}{\lambda}},\lambda\,z,{\lambda}^
{4}R) 
$$
and hence the base $\mathbf{B}$ is a Zariski-open subset of the \it weighted projective 3-space \rm with quasi-homogeneous coordinates $(p:q:r:s)$. The 2-element group acts by the involution
$$
  (p,q,r,s,z,R)\mapsto(s,r,q,p,\frac{1}{z},\frac{1}{R}).
$$
So all the strata of the critical filtration under consideration lie in the factor of weighted projective 3-space over the involution.\\
\\
The stratum $\mathcal{{HUR}}_{4,0;5}^\square$ is the union of 4 irreducible surfaces;  it turns out that
all of them are \underline{rational} and all the Fried families can be calculated explicitly. It would be interesting to know the bounds of this phenomenon in terms of degree.\\
\\
\bf3.1. Genus 1, degree 3. \rm We analyze this case mostly by the pure classical algebraic geometry.\\
\\
Think about the curves of genus 1 as about plain cubics $\mathbf{E}\subset\mathbf{P}_2(\Bbbk)$ and about the functions of degree 3 on it as about \it projections \rm $x_P:\mathbf{E}\to\mathbf{P}_1(\Bbbk)$ from a point 
$P\in\mathbf{P}_2(\Bbbk)\setminus\mathbf{E}$. The best way for such an interpretation of such a  function $x_P\in\Bbbk(\mathbf{E})_3$ is to realize $\mathbf{P}_1(\Bbbk)$ as the \it set of lines \rm in 
$\mathbf{P}_2(\Bbbk)$ through $P$ and to define 
$$
x_P:\mathbf{E}\to\mathbf{P}_1(\Bbbk):Q\mapsto\overline{PQ},
$$ 
where $\overline{PQ}$ denotes the line
through $P$ and $Q$.\\
\\
The dimension of the set of pairs  $(\mathbf{E}\subset\mathbf{P}_2,P\in\mathbf{P}_2)$ is $9+2=11$, but we consider such pairs up to projective equivalence, and 
$$
\dim\frac{\{(\mathbf{E},P)\}}{\mathsf{PGL_3}}=11-8=3,
$$
consistently with the isomorphism
$$
\mathcal{{HUR}}_{3,1}^\square\cong\frac{\{(\mathbf{E},P)\}}{\mathsf{PGL_3}}.
$$
Now fix $\mathbf{E}$ and imagine $P$ moving in $\mathbf{P}_2\setminus\mathbf{E}$. The critical values of $x_P$ correspond to \it tangents \rm from $P$ to $\mathbf{E}$. For a \it generic \rm $P$ there are six of them -- classically this 6 is called the \it class \rm of a cubic curve. Hence we have a critical filtration
$$
\mathcal{{HUR}}_{3,1}^\square=\mathcal{{HUR}}_{3,1;6}^\square\supset
\mathcal{{HUR}}_{3,1;5}^\square\supset\mathcal{{HUR}}_{3,1;4}^\square\supset
\mathcal{{HUR}}_{3,1;3}^\square.
$$
It is geometrically clear that $\mathcal{{HUR}}_{3,1;5}^\square$ consists of the pairs $(\mathbf{E},P)$, in which $P$ lies on an \it inflexion line \rm of $\mathbf{E}$ (two tangent points collide) and that
$\mathcal{{HUR}}_{3,1;4}^\square$ consists of the pairs $(\mathbf{E},P)$, in which $P$ lies on an 
\it intersection \rm of two  inflexion lines of $\mathbf{E}$ (two pairs of tangent points collide).\\
\\
Next suppose that $\text{char}(\Bbbk)\ne2,3$ and use standard equations of the cubic curves. Let (traditionally)  one of the inflection lines to be the ''infinite'' one and the curve $\mathbf{E}$ to be defined by the equation 
$$
y^2=\text{cubic polynomial in } x.
$$
In order to describe $\mathcal{{HUR}}_{3,1;5}^\square$ choose the inflexion line containing  $P$ to be the infinite one; then the pencil of lines through $P$ is a set of parallel lines in the affine plane, and any function$(x,y)\mapsto ax+by+c$ with $b\ne0$ is a function on  $\mathbf{E}$ of degree 3 with 5 critical values.\\
\\
To get from $\mathcal{{HUR}}_{3,1;5}^\square$ to $\mathcal{{HUR}}_{3,1;4}^\square$ we should take the above pencil containing an inflexion line. In the convenient coordinates  the equation looks like
$$
y^2=(1+kx)^2+\frac{x^3}{27},
$$
where the inflexion line $y=1+kx$ is immediately seen. We've got a nice Fried family over
$$
\mathbf{B}=\{k\in\Bbbk\mid 4k^3\ne1\}
$$
with the Fried function 
$$
\Phi:=y-kx-1,
$$
having critical values $\{\infty,0,16\,{k}^{3}-2\pm8\,\sqrt {4\,{k}^{6}-{k}^{3}}$\}. Applying elliptic trick, we get
$$
\boxed{\beta_{{{\text{bas}}}}={\frac { \left( 256\,{k}^{6}-64\,{k}^{3}+1
 \right) ^{3}}{{k}^{3} \left( 4\,{k}^{3}-1 \right) }}}
$$
The dependence  of everything only on $k^3$, not just on $k$, follows from the possibility of substitutions $(k,x)\leftarrow(\rho^{-1}k,\rho x)$ for $\rho\in\sqrt[3]{1}$.\\
\\
Finally, to get from $\mathcal{{HUR}}_{3,1;4}^\square$ to $\mathcal{{HUR}}_{3,1;3}^\square$ we should find such a plane cubic that two of its inflexion lines meet on the third one. The curve of our family over $k=0$,
$$
y^2=1+\frac{x^3}{27},
$$
enjoys this property: the inflexion lines $y=\pm1$ meet in the infinite line. We've met this Belyi pair in the proof of the proposition \bf 1.3.3 \rm (rescaled by $x\leftarrow-3x$).\\
\\
\bf3.2. Genus 2, degree 5. \rm  In the same proof we have seen an ''easy'' Belyi function of degree 5 (the minimal possible) on a curve of genus 2. Brian Birch in $\bf{[Birch1994] }$ has presented a ''hard'' one, namely\footnote{our coordinates are different}
$$
\beta= \left( {z}^{2}+z+1 \right) w+{z}^{5}+\frac{5}{2}\,{z}^{4}+5\,{z}^{3}+5
\,{z}^{2}+\frac{5}{2}\,z+1
$$
on the curve
$$
{w}^{2}={z}^{6}+3\,{z}^{5}+{\frac {29\,}{4}}{z}^{4}+\frac{19}{2}{z}^{3}+{
\frac {29\,}{4}}{z}^{2}+3\,z+1.
$$
The corresponding point in $\mathcal{{HUR}}_{5,2;3}^\square$ lies on two Fried curves in 
$\mathcal{{HUR}}_{5,2;4}^\square$, one of which is defined by (the quarter of)
$$
\Phi=4\,{z}^{5}- \left( 15\,p+20 \right) {z}^{4}+ \left( 10\,{p}^{2}+
20\,p+20 \right) {z}^{3}+
$$
$$+\left( 10\,{p}^{2}+20\,p+
20 \right) {z}^{2}
- \left( 15\,p+20 \right) z$$$$+ \left( 4\,{z}^{2}- \left( 12+8\,p
 \right) z+4 \right) w+4
 $$
on the curves
$$
{w}^{2}={z}^{6}+ \left(\frac{7}{2}\,p-4 \right) {z}^{5}+ \left( {\frac {17\,
}{16}}{p}^{2}-\frac{3}{2}\,p \right) {z}^{4}+$$$$+ \left( -\frac{1}{2}\,{p}^{3}+{\frac {7\,
}{8}}{p}^{2}+4\,p+10 \right) {z}^{3}+ \left( {\frac {17\,}{16}}{p}^{2}
-\frac{3}{2}\,p \right) {z}^{2}+$$$$- \left( \frac{7}{2}\,p+4 \right) z+1.
$$
The elliptic trick provides
$$
\boxed{
\beta_{{{\it bas}}}={\frac {\mathcal{N}^3}{{p}^{4} \left( 2\,p+5 \right) ^{6} \left( {
p}^{2}-4\,p-16 \right)  \left( p-2 \right) ^{3} \left( p+2 \right) ^{5
}}}
}
$$
with
$$
\mathcal{N}=   {p}^{10}-40\,{p}^{8}-80\,{p}^{7}+
320\,{p}^{6}+1088\,{p}^{5}+$$$$+320\,{p}^{4}-1440\,{p}^{3}+720\,{p}^{2}+
5120\,p+4096. 
$$
Birch point corresponds to $p=-2$.\\
\\
It would be interesting to analyze this case geometrically (like the cubics in \bf3.1\rm), realizing genus-two curves as space quintics.\\
\\
\bf3.3. Dreaming about applications. \rm The last example is included in order to show the possible relations of our constructions with other domains of mathematics. \\
\\
Belyi pairs and their deformations for more than a decade appear in several domains of the theory of ordinary differential equations -- algebraic solutions of Painlev\'e VI equations, transformation of hypergeometric-to-Heun equations, etc.; see  $\bf{[SekiVidu2017] }$ and references therein. \\
\\
One of the families this origin ($\bf{[SekiVidu2017] }$) produces the Fried map on the trivial family
$$
\Phi={\frac { \left( s{x}^{3}+15\,{x}^{2}+20\,x+8 \right) ^{2}}{64\,
 \left( x+1 \right) ^{5}}}.
$$
Here $s\in\mathbf{B}=\mathbf{P}_1\setminus{3}$ and the fibers are also rational with the affine coordinate $x$.\\
\\
One calculates
$$
\boxed{
\beta_{{{\mathrm{bas}}}}=-{\frac {s \left( 108\,{s}^{2}-700\,s+1125
 \right) ^{2}}{50000\, \left( s-3 \right) ^{3}}}
 }
$$
 By the way, the factorization of coefficients $$1125=3^2\cdot5^3$$$$729=3^6$$and$$ 2160=2^4\cdot3^3\cdot5$$ points out three special primes. It would be interesting to consider them from the viewpoint of well-developed theory of $p$-adic differential equations.\\
\normalsize
\\
One checks
$$
{\mathrm{div}} \left( \beta_{{{\mathrm{bas}}}} \right) =A_{{1}}+2\,A_{{2}}+2\,A_
{{3}}-3\,C_{{1}}-2\,C_{{2}}
$$
and
$${\mathrm{div}} \left( \beta_{{{\mathrm{bas}}}}-1 \right) =4\,B_{{1}}+B_{{2}}-3\,
C_{{1}}-2\,C_{{2}},
$$
where
$$
s(A_1)=0, s(A_{{2},3}) ={\frac {175}{54}}\pm{\frac {5\,\sqrt {10}}{54},}
$$
$$
s(B_1)=\frac{5}{2}, s(B_2)=\frac{80}{27},
$$
$$
s(C_1)=3, s(C_2)=\infty.
$$
The author hopes that the structures on the base of the families 
 like this one can hope to understand will be useful in understanding the algebraic nature of some differential equations.\\
\begin{center}
\bf4. Conclusion
\end{center}
We indicate briefly several directions of the further research related to the objects discussed in this paper.\\
\\
\bf Combinatorics. \rm Counting Belyi pairs over $\mathbb{C}$ is equivalent to counting dessins d'enfants, and lots of recent papers are devoted to it; they contain precise results (mostly in terms of generating functions) and demonstrate diverse relations with various domains of mathematics and physics. The corresponding counting problems related to critical filtrations should be formulated and solved.\\
\\
\bf Arithmetic. \rm Even in the few explicit examples of the paper the special behavior of certain small primes is obvious. The theoretical complete theoretical explanation of this phenomenon is needed, both in terms of the lists of \it bad \rm primes and the quantitative measures of \it how bad \rm they are -- e.g., it would be nice  to predict combinatorially the powers of primes in the denominators of norms of $\bf{j}$-invariants. The distribution of Galois-orbits of Belyi pairs along the components of the Fried curves should also be studied. \\
\\
\bf Geometry. \rm The components of strata of the critical filtrations should be projected to moduli spaces and studied in terms of  their structures. This study includes the intersections, the behavior of the strata on the boundaries (it demands the extension of the construction of this paper to  stable curves), homology of the strata, relations with level structures, etc.\\
\\
\bf Teichmuller theory. \rm Over $\mathbb{C}$ there is a distinguished homotopic equivalence between two \it close \rm curves on the same component of a Fried curves; therefore, there is a well-defined extremal quasiconformal map between them and hence the distinguished quadratic differentials (defined up to a real positive scalar) on both of them. Are there the cases in which these quadratic differentials have an arithmetic nature?\\
\\
The author together with his students hopes to investigate some of the listed problems and would be happy to involve any interested mathematician.\\
\normalsize
\\
\bf References\rm\\
\\
$\bf{[BaileyFried2002] }$ P. Bailey and M.D. Fried, \it Hurwitz monodromy, spin separation and higher levels of a Modular tower\rm, Arithmetic fundamental groups and noncommutative algebra. PSPUM vol. 70 of AMS(2002), 79-220.
\\
$\bf{[Belyi1980] }$ G.V. Belyi, \it Galois extensions of a maximal cyclotomic fields. \rm Mathematics of the USSR Izvestiya 14, no. 2(1980), 247-256.
\\
$\bf{[Belyi2002] }$ G.V. Belyi, \it A new proof of the three-point theorem. \rm Mathem. Sb. 193, no.3, 21-24.
\\
$\bf{[Birch1994] }$ B. Birch, \it Noncongruence Subgroups, Covers and Drawings\rm, in \it The Grothendieck Theory of Dessins d'Enfants\rm, LMS Lecture Note Series 200, ed. by L. Schneps, Cambridge University Press, 1994, pp. 79--113.
\\
$\bf{[ Bouscaren1998] }$ Elisabeth Bouscaren, ed. \it Model theory and algebraic geometry. \rm . Lecture Notes in
Mathematics, vol. 1696. 
Berlin: Springer-Verlag, 1998.
\\
$\bf{[Clebsch1872] }$ A. Clebsch, \it Zur Theorie der Riemann’schen Flachen. \rm Math. Ann. 6 (1872),
p. 216–230.
\\
$\bf{[Deopurkar2014] }$ Anand Deopurkar, \it Compactifications of Hurwitz spaces. \rm
International Mathematical Research Notices, 2014(14):3863–3911, 2013.
\\
$\bf{[DremShab2017] }$ V. Dremov, G. Shabat, \it On the isomorphisms of bases of Hurwitz families\rm. In preparation.
\\
$\bf{[Fried1977] }$ M. Fried, \it Fields of definition of function fields and Hurwitz families and groups as
Galois groups\rm. Communications in Algebra 5 (1977), 17–82.
\\
$\bf{[Fried1990] }$ M. Fried,  \it Arithmetic of 3 and 4 branch point covers: a bridge provided by noncongruence subgroups of SL2(Z)\rm. Progress in Math. Birkhauser 81 (1990), 77–117. 
 \\
$\bf{[GriffHarr1978]}$ Phillip Griffiths and Joseph Harris, \it Principles of Algebraic Geometry. \rm Wiley Interscience, 1978.
\\
$\bf{[Grothendieck1984] }$ Alexander Grothendieck,  \it Esquisse d'un Programme, \rm (1984 manuscript),  published in Schneps and Lochak (1997, I), pp.5-48; English transl., ibid., pp. 243-283.
\\
$\bf{[Hess2002] }$ F. Hess, \it Computing Riemann-Roch spaces in algebraic function fields and related topics\rm. J. Symboiic Computation, vol. 33, issue 4, April 2002, pp. 425-445.
\\
$\bf{[Hurwitz1891] }$ A. Hurwitz, \it Uber Riemann’sche Flachen mit gegebenen Verzweigungspunkten. \rm
Math. Ann. 39 (1891), p. 1–61.
\\
$\bf{[LandoZvonkin2004] }$ Sergei  Lando and  Alexander Zvonkin,  \it Graphs on Surfaces and Their Applications. 
\rm Springer-Verlag, 2004.
\\
$\bf{[Robinson1963] }$ Robinson, Abraham, \it Introduction to model theory and to the metamathematics of algebra. \rm Amsterdam: North-Holland,  1963.
\\
$\bf{[RomagnyWewers2006] }$ Matthieu Romagny and Stefan Wewers. \it Hurwitz spaces. \rm In Groupes de Galois arithm\'etiques
et diff\'erentiels, volume 13 of S\'emin. Congr., pages 313–341. Soc. Math. France, Paris, 2006.
\\
$\bf{[SekiVidu2017] }$ Jiro Sekiguchi and Raimundus Vidunas, \it Differential relations for almost Belyi maps. \rm To appear (THIS VOLUME???)
\\
$\bf{[Shabat2008] }$ George B. Shabat, \it Visualizing Algebraic Curves:
from Riemann to Grothendieck. \rm Journal of Siberian Federal University. Mathematics and Physics 1 (2008) 42-51.\\
$\bf{[ShabatVoevodsky1990] }$ G.B. Shabat, V.A. Voevodsky, \it
Drawing curves over number fields. \rm
P. Cartier, L. Illusie, N. Katz, G. Laumon, Y. Manin, K. Ribet (Eds.), The Grothendieck Festschrift (5th ed.), Vol. 3, Birkhauser, Basel (1990), pp. 199–227.\\
$\bf{[Xiao2009] }$ Xiao R, Rayson P, McEnery A,  \it A frequency dictionary of mandarin Chinese: Core vocabulary  for learners. \rm London: Routledge; 2009.
\end{document}